\documentclass[12pt]{article}
\usepackage{amsthm}
\usepackage{amssymb,latexsym,amsmath}

\oddsidemargin=\evensidemargin \footskip=36pt
\DeclareMathOperator{\Mat}{Mat}

\DeclareMathOperator{\Irr}{Irr}
\def\css{\begin{cases}}
\def\ecss{\end{cases}}

\def\C{{\mathbb C}}

\def\nmrt{\begin{enumerate}}
\def\enmrt{\end{enumerate}}
\def\tm#1{\item[{\rm (#1)}]}

\def\CC{{\cal C}}

\def\bull{\vrule height 1.2ex width 1ex depth -.1ex }

\makeatletter
\renewcommand{\subsection}{\@startsection{subsection}{2}{0mm}{-2mm}{-2mm}
{\bf\normalsize}}

\makeatother
\newtheorem{formula}{}[section]

\newtheorem{definition}[formula]{Definition}
\newtheorem{corollary}[formula]{Corollary}
\newtheorem{remark}[formula]{Remark}
\newtheorem{lemma}[formula]{Lemma}
\newtheorem{theorem}[formula]{Theorem}

\newtheorem{example}[formula]{Example}

\begin{document}

\title{ }
\author{ J. Bagherian\\ \\ \\ }
\smallskip
\makeatletter
\medskip
\makeatother

\title{Burnside-Brauer Theorem and Character Products in Table Algebras}
\author{
J. Bagherian  \\
Department of Mathematics, University of Isfahan,\\
P.O. Box: 81746-73441, Isfahan, Iran.\\
bagherian@sci.ui.ac.ir \\\\
 A. Rahnamai Barghi
\footnote{Corresponding author: rahnama@kntu.ac.ir}\\ Department
of Mathematics, Faculty of Science
\\ K. N. Toosi University of Technology \\ P.O. Box: 16315-1618,
Tehran, Iran.}

\maketitle
\begin{abstract}
\noindent
In this paper, we first show that the irreducible
characters of a quotient table algebra modulo a normal closed
subset can be viewed as the irreducible characters of the table
algebra itself. Furthermore, we define the character products for
table algebras and give a condition in which  the products of two
characters are characters. Thereafter, as a main result we state
and prove the Burnside-Brauer Theorem on finite groups for table
algebras.
\end{abstract}
\smallskip


\section{{Introduction}}
\indent

One of important results in the character theory of finite groups is the
Burnside-Brauer Theorem. This theorem stats that if a finite group
$G$ has a faithful character $\chi$ which takes $k$ values on $G$,
then every irreducible character of $G$ is a constituent of one of
the characters $\chi^i$ for $0\leq i < k$. One of important
results in this paper is to state and prove an analog of the
Burnside-Brauer Theorem for table algebras. Therefore, we deal
with  products of characters in table algebras. We mention that
products of characters in table algebras need not be  characters
in general. In order to provide a condition in which the products of characters of a given table algebra
are characters, we need to observe the relationship  between the
characters of a table algebra and the characters of its quotient.

In section \ref{prelim}, some elementary facts about table algebras are given. Section
\ref{charquotient} deals with the characters of the quotient
table algebras. More precisely, for a given
table algebra $(A,B)$ and a normal closed subset
$C$ of $B$, we show that the set of irreducible complex characters
of $(A/C,B/C)$ can be embedded in the set of irreducible complex
characters of $(A,B)$. This is a generalization of \cite[Section
3]{Hana1} for association schemes to table
algebras.

Another interesting problem on characters of table
algebras is character products. Since table algebras are not Hopf
algebras in general, character products need not be characters. In
\cite{Hana3}, Hanaki defined character products for
association schemes and gave a condition which  implies that character
products are characters. In section \ref{charproduct}, we define
the character products for table algebras and by
using the results in Section \ref{charquotient} we obtain a
condition for which  character products are characters.
Finally, we prove the Burnside-Brauer Theorem for
table algebras which is a well know theorem in the theory of
finite groups.

\section{Preliminaries}\label{prelim}
\indent
Throughout this paper we follow from \cite{Ar2} for the definition
of  {\it non-commutative table algebras} and related notions.
Hence we deal with non-commutative table
algebras as the following:\\
\medskip
A non-commutative table algebra $(A,B)$ is a finite dimensional
algebra $A$ over the complex field $\C$ and a distinguished basis
${B}=\{b_1=1_A,\cdots,b_d\}$ for $A$, where $1_A$ is a unit,  such
that
the following properties hold:\\
\nmrt \tm{I} The structure constants of ${B}$ are nonnegative real
numbers, i.e., for $a,b \in {B}$:
$$
ab = \displaystyle\sum _{c\in {B}}\lambda_{abc}c,\ \ \ \ \
\lambda_{abc}\in \mathbb{R}^+\cup\{0\}.
$$
\tm{II} There is a semilinear involutory anti-automorphism
(denoted by $^{*}$) of $A$ such that ${B}^{*} = {B}$. \tm{III} For
$a,b\in B$ the equality $\lambda_{ab1_A} = \delta_{ab^*}|a|$ holds
where $|a|>0$ and $\delta$ is the Kronecker symbol. \tm{IV} The
mapping $b\rightarrow |b|, b\in B$ is a one-dimensional linear
representation of the algebra $A$ such that $|b| = |b^*|$ for all $b\in B$
which is called the {\it degree
map}. \enmrt

Throughout this paper a {\it table algebra} means a non-commutative table algebra.\\

Let $(A,B)$ be a table algebra. Then \cite[Theorem 3.11]{Ar}
implies that $A$ is semisimple.  The value $|b|$ is called the
{\it degree} of the basis element $b$. From condition (IV) we see
that $|b| = |b^*|$ for all $b\in B$. Therefore, for an arbitrary
element $\sum_{b\in B}x_b b \in A$, we have $|\sum_{b\in
B}x_bb|=\sum_{b\in B}x_b|b|$.

For each $a=\sum_{b\in B}x_bb$, we set $a^*=\sum_{b\in
B}\overline{x_b}b^*$, where $\overline{x_b}$ means the complex
conjugate of $x_b$. For any $x=\sum_{b\in B}x_bb\in A$, denote by
${\rm{Supp}}(x)$ as the set of all basis elements $b\in B$ such
that $x_b\neq 0$. If $E,D\subseteq B$, then we set
$ED=\bigcup_{e\in E,d\in D }{\rm{Supp}}(ed)$.

A nonempty subset $C\subseteq B$ is called a {\it closed subset},
if ${C}^{*}{C} \subseteq {C}$. We denote by $\mathcal{C}(B)$ the
set of all closed subsets of $B$. In addition,
$C\in\mathcal{C}(B)$ is said to be {\it normal} in $B$ if $bC=Cb$
for every $b\in B$, and denote it by $C\unlhd B$.

Let $(A,B)$ be a table algebra with basis ${B}$ and let $C\in
\mathcal{C}(B)$. From \cite[Proposition 4.7]{Ar}, it follows that
$\{{C}b{C} \mid \ b\in {B}\}$ is a partition of ${B}$. A subset
${C}b{C}$ is called a {\it $C$-double coset} or {\it double coset}
with respect to the closed subset ${C}$. Let
$$
b/C := |C^+|^{-1}({C}b{C})^{+} = |C^+|^{-1}\displaystyle\sum_{x\in
{C}b C}x
$$
where $C^+ = \sum_{c\in C}c$ and  $|C^+|=\sum_{c\in C}|c|$. Then
the following theorem is an immediate consequence of \cite[Theorem
4.9]{Ar}:

\begin{theorem}
\label{Arad} Let $(A,{B})$ be a table algebra and let $C\in
\CC(B)$. Suppose that $\{b_1=1_A,\ldots ,b_k\}$ be a complete set
of representatives of $C$-double cosets. Then the vector space
spanned by the elements $b_i/{C}, 1\leq i\leq k$, is a table
algebra ( which is denoted by $A/{C}$) with a distinguished basis
${B}/{C} = \{b_i/{C} \mid \ 1\leq i \leq k \}.$ The structure
constants of this algebra are given by the following formula:
$$\gamma_{ijk}=
|C^+|^{-1}\sum_{r\in Cb_iC,s\in Cb_jC} \lambda_{rst} $$ where
$t\in C b_k C$ is an arbitrary element.\hfill \bull
\end{theorem}

The table algebra $(A/C,B/C)$ is called the {\it quotient table
algebra} of $(A,B)$ modulo $C$.

We refer the reader to \cite{Zi} for the background of association
schemes.

\section{Embedding of $\Irr(A/C)$ into $\Irr(A)$}\label{charquotient}
In this section we show that there is an embedding from
$\Irr(A/C)$ into $\Irr(A)$ where $C\unlhd B$. We mention that such
an embedding is given for association schemes in \cite{Hana1}.
\\\\
Let $(A,B)$ be a table algebra and $C\in \CC(B)$. Set
$e=|C^+|^{-1}C^+$. Then $e$ is an idempotent for the table algebra
$A$ and the subalgebra $eAe$ is equal to the quotient table
algebra $(A/C,B/C)$ modulo $C$, see \cite{Ar}.
\begin{lemma}
\label{idm} Let $C \in \mathcal{C}(B)$. Then $C\unlhd B$ if and
only if $e=|C^+|^{-1}C^+$ is a central idempotent of $A$.
\end{lemma}
\begin{proof}
Let $C$ be a closed subset of $B$. Clearly $e^2=e$. We first
assume that $C\unlhd B$. In order to prove that $e$ is central, it
is enough to show that for every $b\in B$, $bC^+=C^+b$. From the
normality of $C$ it follows that $(bC)^+=(C b)^+$. From
\cite[Proposition 4.8(ii)]{Ar}, we have $(C b)^+=\alpha^{-1} C^+
b$ and $(bC)^+=\beta^{-1}b C^+$, for some suitable $\alpha, \beta
\in \mathbb{R}^*$. Therefore $\alpha^{-1} C^+ b=\beta^{-1}b C^+$.
But
 $|C^+b|=|C^+||b|=|b||C^+|=|bC^+|$ and so $\alpha=\beta$
which forces $bC^+=C^+b$, as desired.

Conversely let $e$ be a central idempotent. Then
$$
Cb={\rm{Supp}}(C^+b)={\rm{Supp}}(bC^+)=bC
$$
which means that $C\unlhd B$, and we are done.\hfill\bull
\end{proof}

\begin{lemma}
\label{pi} Let $C\unlhd B$. Then the map $\pi: A\rightarrow A/C$
defined by $\pi(b)=\alpha_{b}(b/C)$ where $\alpha_b \in
\mathbb{R}$ such that $C^+ b=\alpha_{b}(C b)^+$, is an algebra
homomorphism.
\end{lemma}
\begin{proof}
Set $e=|C^+|^{-1}C^+$. Then from Lemma \ref{idm}, $e$ is a central
idempotent. We define the map $\pi:A\rightarrow eAe$ where
$\pi(b)=ebe$. It is easily seen that $\pi$ is an algebra
homomorphism. From \cite[Proposition 4.8(ii)]{Ar} there is
$\alpha_{b}\in \mathbb{R}^*$ such that $C^+ b=\alpha_{b}(C b)^+$.
From this fact along with the obvious equality $ebe =
|C^+|^{-1}C^+ b$ we deduce that $ebe=\alpha_{b}|C^+|^{-1}(C
b)^+=\alpha _b (b/C)$. Thus $\pi:A\rightarrow A/C$ where
$\pi(b)=\alpha_{b} (b/C)$ is an algebra homomorphism and we are
done.\hfill\bull
\end{proof}

\begin{theorem}
\label{embbed} Let $C\unlhd B$ and let $\psi:A/C \rightarrow
\rm{Mat_{s}(\C)}$ be a representation of $A/C$. Then
$\overline{\psi}:A\rightarrow \rm{Mat_{s}(\C)}$ defined by
$\overline{\psi}(b)=\alpha_{b}\psi(b/C)$, where $C^+ b=\alpha_{b}(C
b)^+$, is a representation of $A$.
\end{theorem}
\begin{proof}
Put $\overline{\psi}=\psi\circ\pi$, where $\pi$ is defined in
Lemma \ref{pi}. Then $\overline{\psi}$ is an algebra homomorphism
and $\overline{\psi}(b)=\psi\circ\pi(b)
=\psi(\alpha_{b}b/C)=\alpha_{b}\psi(b/C)$, as desired.\hfill\bull
\end{proof}
\vspace{.5cm}
\begin{remark}\label{231207r}
\nmrt \tm{1} By Theorem \ref{embbed} we may embed $\Irr(A/C)$ into
$\Irr(A)$, indeed if $\chi$ and $\psi$ are  distinct characters of
$A/C$, then $\chi(b/C) \neq \psi(b/C)$, for some $b/C \in B/C$ and
so $\overline{\chi}(b) \neq \overline{\psi}(b)$, where $\overline{\chi}$
is  a representation of $A$ defined by $\overline{\chi}(b) = \alpha_b\chi(b/C)$.
\tm{2} In Theorem
\ref{embbed}, the correspondence preserves the irreducibility of
representations. Furthermore, if $D$ is an irreducible
representation of $A$ such that rank $D(e)\neq 0$, then
$ D = \overline{E}$ where $E$ is appropriately explained. \enmrt
\end{remark}
In Theorem \ref{embbed}, the value of $\alpha_b$ is not be given
precisely. But if we consider another property on $C$ which is
stronger than normality, namely {\it strongly normal closed
subset}, then it is possible to give the value of $\alpha_b$ in
precise form. The rest of this section deals with the quotient of
table algebras modulo $C$, where $C$ is a strongly normal closed
subset.

\begin{definition}
A closed subset $C$ of $B$ is said to be {\it strongly normal} and
denoted by $C\unlhd' B$,  if for each $b\in B$
$$
b^* C b\subseteq C.
$$
\end{definition}

In the following we show that a strongly normal closed subset is a
normal closed subset.
\begin{lemma}
Every strongly normal closed subset is a normal closed subset.
\end{lemma}
\begin{proof}
Let $C \unlhd' B$. Then for $b\in B$ we have $b^*C b\subseteq C$,
and so $bC b^*C b\subseteq bC$ which implies that $C b\subseteq
bC$. On the other hand, from $bC b^*\subseteq C$ it follows that $bC b^*C
b\subseteq C b$, and so $bC \subseteq C b$. Thus $C b = bC $, as
desired.\hfill\bull
\end{proof}
\vspace{0.4cm} The following example shows that the converse of
the above lemma is not true, i.e., a normal closed subset is not
necessarily strongly normal closed subset.

\begin{example}
Let $q \geq 2$ and $B = \{r_0,r_1,\ldots ,r_{q+1} \}$ be a basis
for a complex vector space $A$ of dimension $q+2$. We define
multiplication
\[r_ir_j=\left\{\begin{array}{ll} (q-1)r_0+(q-2)r_i, & \text{if}~~~ i=j\neq 0\\
\displaystyle\sum_{k\neq {0,i,j}}r_k, & \text{if}~~~{i\neq
j}\end{array}\right.\]
for all $i,j$ with $1\leq i,j \leq q+1$, and $r_ir_0 = r_i$ for all $i$.
This extends to a multiplication in $A$
which is commutative and association with unit element $r_0 =
1_A$. A direct computation shows that $(A,B)$ is a table algebra
where $r_i^* = r_i$ for all $i$; and $|r_i| = q-1$ for $i>0$.
Clearly the set $\{r_0,r_i\}$ for every $i\neq 0$ is a normal
closed subset but it is not a strongly normal closed subset. In
fact, the construction of this table algebra is given in
\cite{PR}.
\end{example}

\begin{theorem}\label{gr}
Let $b\in B$ and $C \in \mathcal{C}(B)$. Then $|b/C|=1$ if and
only if $b^*C b \subseteq C$.
\end{theorem}
\begin{proof}
Let $T=\{b_1=1_A, b_2,\ldots,b_t\}$ be a complete set of
representatives of $C$-double cosets and let $b=b_i$ for some $i,~
1\leq i \leq t$.

Suppose that $|b/ C|=1$ and  let $d\in b^*Cb$. Since $1_A\in
C$, we have $d \in (C b^*C )(C bC)$. Then there exists $r \in
Cb^*C$ and $s\in C bC$ such that $\lambda_{rsd}\neq 0$. As
$B=\bigcup_{j = 1}^t Cb_jC$, we may assume that $d\in C b_k C$
for some $k,~~ 1\leq k \leq t$. So $\gamma_{i^*ik}\neq 0$ where
$$
\gamma_{i^*ik}=|C^+|^{-1}\displaystyle\sum_{r\in C {b_i}^*C,s\in C
b_iC} \lambda_{rsd}.
$$
But from the assumption we conclude that
$(b/C)(b^*/C)=\{1_A/C\}$. Hence $\gamma_{i^*i1}\neq 0$ and so
that $k=1$. Thus $d\in C$ and so $b^*C b \subseteq C$.

Conversely, let $b^*C b \subseteq C$. Then $(Cb^*C) (C
b C)\subseteq C$ and so
\begin{eqnarray}\label{181207}
\gamma_{i^*i1}\leqslant |C^+|^{-1}\displaystyle\sum_{r, s \in C}
\lambda_{rsd}
\end{eqnarray}
for  $d\in C$. Now from (\ref{181207}) and the equalities
$\displaystyle\sum_{r,s\in C}
\lambda_{rsd}=\displaystyle\sum_{r\in C} |r|=|C^+|$ (see
\cite[Proposition 2.3(i)]{Ar}) we deduce that $\gamma_{i^*i1}=1$.
Thus $|b/C|=1$, as desired. \hfill\bull
\end{proof}
\begin{corollary}\label{group}
Let $(A,B)$ be a table algebra and $C\in \mathcal{C}(B)$. Then
$(A/C,B/C)$ is a group algebra if and only if $C\unlhd' B$.
\end{corollary}
\begin{proof}
This follows immediately from Theorem \ref{gr}. \hfill \bull
\end{proof}

\begin{theorem}
\label{mainembed} Let $C\unlhd' B$ and $\psi$ be a representation
of $A/C$. Then the mapping $\overline{\psi}: b \mapsto
|b|\psi(b/C)$ for every  $b\in B$, is a representation for $A$.
\end{theorem}
\begin{proof}
Using the notation in the proof of Theorem \ref{embbed} it
suffices to show that  $\alpha_b=|b|$. Since $C\unlhd' B$, it
follows that $bb^*C\subseteq C$. Then ${\rm{Supp}}(bb^*)\subseteq
C$ and from \cite[Proposition 4.8(iv)]{Ar} the equality $|(C
b)^+|=|C^+|$ follows. On the other hand, from \cite[Proposition
4.8(ii)]{Ar} there exists $\alpha_b\in \mathbb{R}^*$ such that
$|C^+|| b|=\alpha_{b}|(C b)^+|$. Thus $|C^+||b|=\alpha_{b}|C^+|$
and so $\alpha_{b}=|b|$, as desired. \hfill\bull
\end{proof}

\section{{Character products}}\label{charproduct}
\indent

For an associative algebra A, the tensor product $V\otimes W$ of
two $A$-modules $V$ and $W$ is a vector space, but not necessarily
an $A$-module. In order to make an $A$-module on $V\otimes W$,
there must be a linear binary operation $\Delta: A\rightarrow
A\otimes A$ which is also an algebra homomorphism. This is an
important property for the algebra $A$ becomes a Hopf algebra. For
instance, in group theory the tensor products of two $G$-modules
$V$ and $W$ gives us a module, indeed the group algebra $\C G$ is
a Hopf algebra with $\Delta: g\rightarrow g\otimes g.$ So if
$\chi$ and $\psi$ are afforded by two $G$-modules, then their tensor
product affords the character $\chi\psi(g):=\chi(g)\psi(g)$ which
is called the character product of $\chi $ and $\psi$.

In general, a table algebra $(A,B)$ is not a Hopf algebra and so
it is not generally possible to define the structure of an
$A$-module on $V\otimes W$. In \cite{Do} Doi introduced a
generalization of Hopf algebras and defined a binary linear
operation $\Delta:b \rightarrow \frac{1}{|b|}b \otimes b$, $b\in
B$. By considering this binary linear operation, we define the
character products of $\chi$ and $\psi$ by:
\begin{eqnarray}\label{271207d}
\chi\psi(b):=\frac{1}{|b|}\chi(b)\psi(b),~~ b\in B.
\end{eqnarray}

Since $\Delta$ is not necessarily an algebra homomorphism, a
character products in a table algebra is not generally a character.
It might be mentioned that this is an analog of association
schemes which has
already done by Hanaki in \cite{Hana3}.\\

Through this section we assume that $(A,B)$ is a table algebra
with a strongly normal closed subset $C$ and $e=|C^+|^{-1}C^+$.\\

Let $V$ and $W$ be $A/C$-module and $A$-module, respectively. We
define a multiplication of $A$ on $V\otimes W$ as the following:
\begin{eqnarray}\label{strmod}
b(v\otimes w):=(b/C)v\otimes bw , ~~ v\in V, ~w\in W, ~ b\in B.
\end{eqnarray}

\begin{lemma}\label{221207}
Let $V$ be  an irreducible $A$-module with $\dim_\C(eV)\neq 0$ and
let $W$ be an $A$-module. Then $V\otimes W$ is an $A$-module given
by the multiplication in (\ref{strmod}).

\end{lemma}
\begin{proof}
We first claim that $\mu:A\rightarrow A/C\otimes A$ by $\mu(b)=b/C
\otimes b$, for any $b\in B$
 is an algebra homomorphism.
Let $b,c\in B$ be given. Then
\begin{eqnarray}\label{290108}
\mu(bc)=\mu(\displaystyle\sum_{d\in
B}\lambda_{bcd}d)=\displaystyle\sum_{d\in B}\lambda_{bcd}\mu(d)=
\displaystyle\sum_{d\in B}\lambda_{bcd}d/C\otimes d.
\end{eqnarray}
On the other hand,
$$
\mu(b)\mu(c)=(b/C \otimes b)(c/C \otimes c)
=\displaystyle\sum_{d\in B}\lambda_{bcd}(b/C)(c/C)\otimes d.
$$
Now in the above equality, if $\lambda_{bcd}\neq 0$ then
$\gamma_{b/C,c/C,d/C}\neq 0$ and so $(b/C)(c/C)=d/C$, indeed by
Corollary \ref{group} $B/C$ is a group. Thus
$\mu(b)\mu(c)=\sum_{d\in B}\lambda_{bcd}d/C\otimes d$ which is
equal to $\mu(bc)$ by (\ref{290108}), and the claim is proved.
Since $\dim_\C(eV)\neq 0$, then from \cite[Corollary 3.6]{BR}, $V$
is an $A/C$-module. This implies that $b(v\otimes
w)=\mu(b)(v\otimes w)$, and the proof is complete. \hfill\bull
\end{proof}

\begin{lemma}\label{product1}
Let $\chi$  be an irreducible  character afforded by $A$-module
$V$ such that $\dim_\C(eV)\neq 0$. Then for every $A$-module $W$
the tensor product $ V\otimes W$ is an $A$-module which affords
the character $\frac{1}{|b|}\chi(b)\psi(b)$.
\end{lemma}
\begin{proof}
We first note that under our assumptions $V\otimes W$ is an
$A$-module, by Lemma \ref{221207}. Since  $\dim_\C(eV)\neq 0$,
from Corollary \cite[Corollary 3.6]{BR} it follows that the module
$V$ can be considered as an $A/C$-module. Now let
$D_1:A/C\rightarrow \Mat_{d_1}(\C)$ and $D_2:A \rightarrow
\Mat_{d_2}(\C)$ be representations of $A/C$ and $A$ corresponding
to $V$ and $W$, respectively. Then $(D_1\otimes D_2)\circ
\mu:A\rightarrow \rm{Mat}_{d_1d_2}(\C) $ is a representation of
$A$ corresponding to $V\otimes W$, where $\mu:A\rightarrow
A/C\otimes A$ is the algebra homomorphism given in the proof of
Lemma \ref{221207}. Now we conclude that the trace of
$((D_1\otimes D_2)\circ \mu)(b) $ is equal to $\chi(b/C)\psi(b)$,
for $b\in B$. But from Theorem \ref{mainembed} and  Remark
\ref{231207r} (2) it follows that $\chi(b/C) =
\frac{1}{|b|}\chi(b)$, and we are done. \hfill\bull
\end{proof}

\begin{theorem}
\label{product} Let $\chi$  be an irreducible character of $A$
such that $\chi(e)\neq0$. Then  $\chi\psi$ is a character of $A$,
where $\psi$ is a character of $A$.
\end{theorem}
\begin{proof}
Since $\chi(e) \neq 0$, it follows from \cite[Corollary 3.5]{BR}  that
$\chi$ is an irreducible character of $A/C$. Thus from  Lemma
\ref{product1} we conclude that the character product $\chi \psi$
is a character, as desired.\hfill \bull
\end{proof}
\begin{theorem}
Let $\chi, \psi\in \Irr(A)$. If $\chi(e)=1$, then $\chi\psi\in
\Irr(A)$.
\end{theorem}

\begin{proof}
Let $V$ and $W$ be two irreducible $A$-modules which afford $\chi$
and $\psi$ respectively. The equality $\chi(e) = 1$ implies that
$\dim _{\C}(V)=\chi(e)=1$. From \cite[Corollary 3.5]{BR} and
Theorem \ref{product} it follows that $\chi$ is a linear character
of $A/C$ and $\chi\psi$ is a character, respectively. Moreover,
$\dim _{\C}(V)=1$ implies that every $A$-submodule of $V\otimes W$
is of the form $V\otimes W'$ where $W'$ is an $A$-submodule of
$W$. Therefore, by irreducibility of $W$ the $A$-module $V\otimes
W$ is irreducible and so $\chi\psi\in \Irr(A)$. This completes the
proof. \hfill \bull
\end{proof}
\vspace{3mm} Let $A$ be a finite dimensional algebra with a basis
$w_1,\ldots,w_r$ over a field $F$. Let $\zeta $ be a
non-degenerate feasible trace on $A$. Then from \cite{Hi}, $\zeta$
induces a {\it dual form} $[\cdot,\cdot]$ on ${\rm Hom}_F(A,F)$ in
which for every $\chi,\varphi\in {\rm Hom}_F(A,F)$ we have
$$
[\chi,\varphi]=\displaystyle \sum_{i
=1}^r\chi(w_i)\varphi(\hat{w_i})
$$
where $\hat{w_1},\ldots ,\hat{w_r}$ is the dual basis defined by
$\zeta(w_i\hat{w_j}) = \delta_{i,j}$.

Now let $(A,B)$ be a table algebra. The linear function $\zeta$ on
$A$ is defined in \cite{BR} by setting $\zeta(b) =
\delta_{b,1_A}|B^+|$, for $b\in B$. Then $\zeta$
is a non-degenerate feasible trace on $A$ and it
follows that the dual form $[\cdot,\cdot]$ on ${\rm Hom}_\C(A,\C)$
is as follows:
\begin{eqnarray}\label{inner}
[\chi,\varphi]=\frac{1}{|B^+|}\displaystyle\sum_{b\in
B}\frac{1}{|b|}\chi(b)\varphi(b^*)
\end{eqnarray}
for every $\chi,\varphi\in {\rm Hom}_\C(A,\C)$.

Let $\chi$ be a character of table algebra $(A,B)$.
The following subset of $B$
$$
K(\chi) = \{ b \in B: \chi(b) = |b|\chi(1) \}
$$
is a close subset of $B$. The proof of this result can be found in
\cite[Theorem 3.2]{Hana0} for the case of association scheme, but the proof
also works for table algebras.

Below we give our main result which proves the Burnside-Brauer
Theorem on finite groups for table algebras.
\begin{theorem}\label{burn}
Let $(A,B)$ be a table algebra. Suppose that $A$ has a character
$\chi$ with $K(\chi) =\{1_A\}$ such that
$\displaystyle\frac{\chi(b)}{|b|}$ takes on exactly $k$ different
values for $b \in B$. If all powers of $\chi$ is a character, then
each irreducible character of $A$ appears as an irreducible
component of one of $\chi^i$, where $0\leq i \leq k-1$ and $\chi^0
= \rho$.
\end{theorem}
\begin{proof}
Let $\alpha_1,\ldots,\alpha_k$ be the distinct values taken by
$\displaystyle\frac{\chi(b)}{|b|} , b\in B$. Define $B_t=\{b\in B
| \chi(b)=|b|\alpha_t\}$. Assume that $\alpha_1=\chi(1)$ so that
$B_1=K(\chi)$. Fix $\psi\in\Irr(A)$ and let $\beta_i
=\displaystyle\sum_{b\in B_i}\psi(b^*)$ for $1\leq i \leq k$.
Since $\chi^2(b)=\frac{1}{|b|}\chi(b)^2$, it follows that
$\chi^j(b)=\frac{1}{|b|^{j-1}}\chi(b)^j$. Hence
 \begin{eqnarray}
[\chi^j,\psi]&=&{\nonumber} \frac{1}{|B^+|}\displaystyle\sum_{b\in
B}\frac{1}{|b|}
 \chi^j(b)\psi(b^*)\\&=&{\nonumber}
 \frac{1}{|B^+|}\displaystyle\sum_{i=1}^{k}\displaystyle\sum_{b\in
B_i}\frac{1}{|b|^{j}}\chi(b)^j\psi(b^*)\\&=&{\nonumber}
\frac{1}{|B^+|}\displaystyle\sum_{i=1}^{k}(\alpha_i)^j\beta_i.
{\nonumber}
\end{eqnarray}
Therefore, if $\psi$ is not a constituent of $\chi^j$ for all
$0\leq j \leq k-1,$ then
 \begin{eqnarray}\label{burnside}
\displaystyle\sum_{i=1}^{k}(\alpha_i)^j\beta_i=0,~~~ j =
0,1,\ldots ,k-1.
\end{eqnarray}
Let $M:=(a_{i,j})$ be a $k\times k$ matrix whose $i$th row and
$j$th column is $(\alpha_i)^j$ and let $X=(\beta_1,\beta_2,\ldots , \beta_k)$.
Therefore (\ref{burnside}) shows that $XM=0$. But the determinant
of $M$ is Vandermonde determinant and is equal to
$\pm\Pi_{i<j}(\alpha_i-\alpha_j)\neq0$. It follows that $X=0$. But
$\beta_1 =\psi(1_A)\neq 0$,
which is a contradiction. This completes the proof of the
theorem.\hfill\bull
\end{proof}
\begin{remark}
By using Theorem \ref{burn} for $(A,B) = (\mathbb{C}G,G)$,
where $G$ is a
finite group,  we get
the Burnside-Brauer Theorem on finite groups.

\end{remark}
\begin{center}
ACKNOWLEDGMENT
\end{center}
The authors thank the referee anonymous for his kind comments and suggestions.


\begin{thebibliography}{9}
\bibitem{Ar2} Z. Arad, E. Fisman, M. Muzychuk, {\em On a Product of Two Elements
in Non-commuatative C-Algebras}, Algebra Colloquium, 5:1, 85-97,
1998.

\bibitem{Ar} Z. Arad, E. Fisman, M. Muzychuk, {\em Generalized Table Algebras, Israel
J. Math.} 114, 29-60, 1999.

\bibitem{BR}
J. Bagherian, A. Rahnamai Barghi, {\em Standard Character Condition for
C-algebras}, arXiv: 0810.5305v1, Submitted to Algebra Colloquium.

\bibitem{Do} Y.Doi, {\em Bi-Frobenius Algebras and Group-Like Algebra, Hopf algebras} in:
Lecture Note in Pure and Appl.Math, Vol. 237, Dekker, New York,
pp.143-155, 2004.

\bibitem {Hana0} A. Hanaki, {\em Characters of association schemes and normal closed subsets}, Graphs Combin. 19, 363-369, 2003.

\bibitem {Hana1} A. Hanaki, {\em Representations of Association Schemes and Their Factor Schemes}, Graphs Combin. 19,
195-201, 2003.

\bibitem{Hana3} A. Hanaki, {\em Character Products of Association Schemes}, J. Algebra 283, 596-603, 2005.

\bibitem{Hi} D. G. Higman, {\em Coherent Algebras}, J. Linear Algebra and it's Applications, 93, 209-239, 1987.

\bibitem{PR} I. Ponomarenko, A. Rahnamai Barghi, {\em On Amorphic C-algebras},
Journal of Mathematical Sciences, Vol. 145, No. 3, 2007.

\bibitem{Zi}P-H. Zieschang, {\em An Algebraic Approach to Association
Schemes}, Lecture Notes in Math., vol. 1628, Springer-Verlag,
 Berlin, 1996.
\end{thebibliography}
\end{document}